\documentclass[twoside]{amsart}
\usepackage{latexsym}
\usepackage{amssymb,amsmath,amsopn}
\usepackage[dvips]{graphicx}   
\usepackage{color,epsfig}      

\newtheorem{thm}{Theorem}[section]
\newtheorem{prop}[thm]{Proposition}

\newtheorem{conj}[thm]{Conjecture}

\theoremstyle{definition}
\newtheorem{defn}[thm]{Definition} 
\newtheorem{rem}[thm]{Remark}

\newtheorem{que}[thm]{Question}
\numberwithin{equation}{section}

\renewcommand{\dim}{d}
\renewcommand{\Re}{\mathbb R}
\renewcommand{\epsilon}{\varepsilon}
\newcommand{\Red}{\Re^\dim}

\newcommand{\dlog}{\dim\log\dim+ \log\log\dim + 5\dim}
\newcommand{\KD}{\mathcal{K}^\dim_a}
\newcommand{\Ha}{H}
\newcommand{\Had}{\Ha_\dim}

\newcommand{\UHad}{\overline{\Ha}_\dim}

\newcommand{\remark}[1]{}

\DeclareMathOperator{\dist}{d_{BM}}
\DeclareMathOperator{\vol}{vol}

\DeclareMathOperator{\card}{card}

\begin{document}
\title[On the Constant of Homothety]{On the Constant of Homothety for Covering a Convex Set with Its Smaller Copies}
\author[M. Nasz\'odi] {M\'arton Nasz\'odi$^*$}

\address{M\'arton Nasz\'odi, Dept. of Math. and Stats., 632 Central Academic Building, University of Alberta, Edmonton, Ab, Canada T6G 2G1}
\email{mnaszodi@math.ualberta.ca}
                                                                               
\thanks{${}^*$ Partially supported by a Postdoctoral Fellowship of the Pacific Institute for the Mathematical Sciences}
\subjclass{52A35, 52A20, 52C17}
\keywords{illumination, Boltyanski--Hadwiger Conjecture, convex sets}

\begin{abstract}
Let $\Had$ denote the smallest integer $n$ such that for every convex body $K$ in $\Red$ there is a
$0<\lambda < 1$ such that $K$ is covered by $n$ translates of $\lambda K$.
In \cite{BMP} the following problem was posed: Is there a $0<\lambda_\dim<1$ depending on $\dim$ only with the property
that every convex body $K$ in $\Red$ is covered by $\Had$ translates of $\lambda_\dim K$?
We prove the affirmative answer to the question and hence show that the Gohberg--Markus--Boltyanski--Hadwiger Conjecture
(according to which $\Had\leq2^d$) holds if, and only if, a formally stronger version of it holds.
\end{abstract}
\maketitle


\section{Definitions and Results}

A \emph{convex body} in $\Red$ is a compact convex set $K$ with non--empty interior. Its volume is denoted by $\vol(K)$.
\begin{defn}
For $\dim\geq 1$ let $\Had$ denote the smallest integer $n$ such that for every convex body $K$ in $\Red$ there is a
$0<\lambda < 1$ such that $K$ is covered by $n$ translates of $\lambda K$. Furthermore, let $\UHad$
denote the smallest integer $m$ such that there is a $0<\lambda_\dim < 1$ with the property that every convex body $K$ 
in $\Red$  is covered by $m$ translates of $\lambda_\dim K$.
\end{defn}
Clearly, $\Had\leq\UHad$. The following question was raised in \cite{BMP} (Problem 6 in Section 3.2): \emph{Is it true that $\Had=\UHad$?}

We answer the question in the affirmative using a simple topological argument.
\begin{thm}\label{thm:univconst}
$\Had=\UHad$.
\end{thm}

The famous conjecture of Gohberg, Markus, Boltyanski and Hadwiger states that $\Had\leq 2^d$ (and only the cube requires $2^d$
smaller positive homothetic copies to be covered). For more information on the conjecture, refer to \cite{B}, \cite{MS} and \cite{Sz}.
In view of Theorem~\ref{thm:univconst}, the conjecture is true if, and only if, the following,
formally stronger conjecture holds:

\begin{conj}[Strong Gohberg--Markus--Boltyanski--Hadwiger Conjecture]
For every $\dim\geq 1$ there is a $0<\lambda_\dim<1$ such that every convex body $K$ in $\Red$ 
is covered by $2^d$ translates of $\lambda_\dim K$.
\end{conj}

In Section~\ref{sec:proof} we prove the Theorem. We note that the proof provides no upper bound on $\lambda_\dim$ in terms of $\dim$.
In Section~\ref{sec:remarks} we show an upper bound on the number of translates of $\lambda K$ 
required to cover $K$, improving a result of Januszewski and Lassak \cite{JL}.

\section{Proof of Theorem~\ref{thm:univconst}}\label{sec:proof}
We define the following function on the set of convex bodies:
\[
\lambda(K):=\inf\{\lambda >0 : K \mbox{ is covered by } \Had \mbox{ translates of } \lambda K\}.
\]

By \cite{R1}, $\Had$ is finite for every $d$, so $\lambda(.)$ is well defined.

\begin{rem}\label{rem:affinv}
Clearly, $\lambda(.)$ is affine invariant; that is, if $T$ is an invertible affine transformation of $\Red$ then
$\lambda(K)=\lambda(TK)$. Moreover, $0<\lambda(K)<1$.
\end{rem}

We recall the definition of the (multiplicative) \emph{Banach--Mazur distance} of two convex bodies $L$ and $K$ in $\Red$:

\begin{equation}\label{eq:BM}
\dist(L,K)=
\end{equation}
\[
\inf\{\lambda > 0 : L-a\subseteq T(K-b) \subseteq \lambda(L-a) \mbox{ for some } a, b\in \Re, T\in GL(\Red)\}
\]

The following proposition states that $\lambda(.)$ is upper semi--continuous. Similar statements have been proved before, cf. Lemma 2. 
in \cite{BS}.
\begin{prop}\label{prop:usemicont}
For every convex body $K$ and $\epsilon>0$ there is a $\delta>0$ with the property that for any convex body $L$,
if $\dist(L,K)<1+\delta$ then $\lambda(L)<\lambda(K)+\epsilon$.
\end{prop}

\begin{proof}
Let $\lambda:=\lambda(K)+\frac{\epsilon}{2}$. Then there is a set $\Lambda\subset\Red$ with $\card \Lambda\leq\Had$ such that
$K\subseteq \Lambda + \lambda K$. Now, let $\delta >0$ be such that 
\begin{equation}\label{eq:delta}
1+\delta< \frac{\lambda+\frac{\epsilon}{2}}{\lambda}
\end{equation}
Assume that $\dist(L,K)<1+\delta$; that is, 
\begin{equation}\label{eq:close}
L-a\subseteq \bar{K}\subseteq (1+\delta)(L-a), 
\end{equation}
where $\bar{K}$ is an affine image (under an invertible affine transformation) of $K$. Clearly, we may assume that $\bar{K}=K$.

It follows that $L-a\subseteq \Lambda+(1+\delta)\lambda (L-a)$, and hence, $\lambda(L)\leq (1+\delta)\lambda < \lambda(K)+\epsilon$.
\end{proof}

Let $\KD$ denote the set of affine equivalence classes of convex bodies in $\Red$ equipped with the topology induced by the
metric $\dist$. In \cite{M} it is shown that $\KD$ is a compact space. (Note that Macbeath uses a different metric on $\KD$ however, that metric
induces the same topology as $\dist$, cf. \cite{Gr}.)

It follows from Remark~\ref{rem:affinv} and Proposition~\ref{prop:usemicont} that $\lambda(.)$ is an upper semi--continuous function on a compact space.
Hence, it attains its maximum, which (by Remark~\ref{rem:affinv}) is less than one. This proves Theorem~\ref{thm:univconst}.

\section{Quantitative Results}\label{sec:remarks}

Januszewski and Lassak \cite{JL} proved that for every $k+l>\dim^\dim$, any convex body $K\subset\Red$ is covered
by $k$ translates of $\lambda K$ and $l$ translates of $-\lambda K$, where $\lambda =1-\frac{1}{(\dim+1)\dim^\dim}$. 
The following argument shows that one may obtain a better estimate on the number of translates of 
$\lambda K$ required to cover $K$, using results of Rogers \cite{R1}, Rogers and Shephard \cite{RS}, and Rogers and Zhong \cite{RZ}.

Let $K, L$ be convex bodies in $\Red$. Let $N(K,L)$ denote the \emph{covering number} of $K$ and $L$; that is, 
the smallest number of translates of $L$ required to cover $K$. In \cite{RZ} it is shown that
\[
N(K,L)\leq\frac{\vol(K-L)}{\vol(L)}\Theta(L),
\]
where $\Theta(L)$ is the covering density of $L$. By \cite{R1}, $\Theta(L)\leq \dlog$ for every convex body $L$
in $\Red$. It follows that for any $0<\lambda<1$ we have
\begin{eqnarray}\label{eq:elso}
N(K,\lambda K)&\leq& \lambda^{-\dim}\frac{\vol(K-K)}{\vol{K}}(\dlog)\nonumber\\
&\leq& \lambda^{-\dim}\binom{2\dim}{\dim} (\dlog)
\end{eqnarray}

The last inequality follows from the Rogers--Shephard Inequality\cite{RS}. Similarily,
\begin{eqnarray}\label{eq:masodik}
N(K,-\lambda K)&\leq& \lambda^{-\dim}\frac{\vol(K+K)}{\vol{K}}(\dlog)\nonumber\\
&=& \lambda^{-\dim}2^\dim(\dlog)
\end{eqnarray}

By substituting $\lambda=\frac{1}{2}$ into (\ref{eq:elso}) and (\ref{eq:masodik}), we obtain the following:
\begin{rem}\label{rem:univexists}
The number of translates of $\frac{1}{2}K$ that cover $K$ is of order not greater than $8^\dim\sqrt{\dim}\log\dim$;
and the number of translates of $-\frac{1}{2}K$ that cover $K$ is of order not greater than $4^\dim\dim\log\dim$.
\end{rem}

\begin{defn}
Let $0<\lambda < 1$, and  $\dim\geq 1$. We denote by $\UHad(\lambda)$ the smallest integer $n$ such that
every convex body $K$ in $\Red$  is covered by $n$ translates of $\lambda K$.
\end{defn}

It follows from Remark~\ref{rem:univexists} that $\UHad(\frac{1}{2})$ is finite for every $d$. A natural strengthening of the question we 
discussed in this note is the following:
\begin{que}
Is there a universal constant $0<\lambda<1$ such that for every dimension $d$, $\Had$ is equal to $\UHad(\lambda)$?
\end{que}

\end{document}